\documentclass{amsart}
\usepackage[utf8]{inputenc}
\usepackage{amsmath}
\usepackage{amsfonts}
\usepackage{amsthm}
\usepackage{amssymb}

\usepackage{graphicx}
\usepackage[dvipsnames]{xcolor}

\usepackage{hyperref}

\usepackage{verbatim}

\title{Schottky groups and maximal representations}
\author{Burelle, Jean-Philippe}
\author{Treib, Nicolaus}
\thanks{The authors acknowledge support from U.S. National Science Foundation grants DMS 1107452, 1107263, 1107367 "RNMS: Geometric Structures and Representation Varieties" (the GEAR Network), and are grateful for the hospitality and the motivating working environment at MSRI during the program "Dynamics on Moduli Spaces of Geometric Structures". The first author was supported by a Postgraduate Scholarship from the Natural Sciences and Engineering Research Council of Canada. The second author was supported by the European Research Council under ERC-Consolidator grant 614733.}

\date{September 2016}

\newtheorem{Prop}{Proposition}
\newtheorem{Cor}{Corollary}
\newtheorem{Lem}{Lemma}
\newtheorem{Thm}{Theorem}

\theoremstyle{definition}
\newtheorem{Ex}{Example}
\newtheorem{Def}{Definition}

\theoremstyle{remark}
\newtheorem*{Rem}{Remark}

\newtheoremstyle{TheoremNum}
    {\topsep}{\topsep}              
    {\itshape}                      
    {}                              
    {\bfseries}                     
    {.}                             
    { }                             
    {\thmname{#1}\thmnote{ \bfseries #3}}
\theoremstyle{TheoremNum}
\newtheorem{thmn}{Theorem}

\newcommand{\GL}[1]{\mathrm{GL}(#1,\bR)}
\newcommand{\Sp}[1]{\mathrm{Sp}(#1,\bR)}
\newcommand{\PSL}[1]{\mathrm{PSL}(#1,\bR)}

\newcommand{\RP}[1]{\mathbb{RP}^{#1}}
\newcommand{\bR}{\mathbb{R}}
\newcommand{\bN}{\mathbb{N}}

\newcommand{\bH}{\mathbb{H}}
\newcommand{\Lag}[1]{\mathrm{Lag}(\bR^{#1})}
\newcommand{\SOtt}{\mathrm{SO}(3,2)}

\newcommand{\M}{\mathsf{M}}
\newcommand{\Sh}{\mathcal{S}}

\newcommand{\Inv}[2]{\sigma_{#1#2}}

\newcommand{\B}[2]{\mathcal{B}_{#1#2}}
\newcommand{\BB}{\mathcal{B}}
\newcommand{\PH}[2]{\mathcal{P}(#1,#2)}
\newcommand{\graph}{\mathrm{graph}}

\newcommand{\cycle}[3]{\overrightarrow{#1#2#3}}

\DeclareMathOperator{\Aut}{Aut}
\DeclareMathOperator{\sgn}{sgn}

\begin{document}
\begin{abstract}
We describe a construction of Schottky type subgroups of automorphism groups of partially cyclically ordered sets. We apply this construction to the Shilov boundary of a Hermitian symmetric space and show that in this setting Schottky subgroups correspond to maximal representations of fundamental groups of surfaces with boundary. As an application, we construct explicit fundamental domains for the action of maximal representations into $\Sp{2n}$ on $\RP{2n-1}$.
\end{abstract}

\maketitle

\section{Introduction}
Let $H\subset \PSL{2}$ be a Fuchsian Schottky group. This means that there exist $\ell_1,\dots,\ell_{2g}$ geodesics in the hyperbolic plane $\bH^2$ with disjoint closures, and a collection of generators $A_1,\dots,A_g\in H$ which identify these geodesics pairwise \cite{OuterCircles}. By a classical ping pong argument, $H$ is free on its generators and the quotient $\bH^2/H$ is a hyperbolic surface with non-empty boundary. The homeomorphism type of the surface depends on the combinatorial type of the pairings. If we allow the geodesics to be asymptotic, we obtain the converse: every complete hyperbolic metric on the interior of a compact surface with boundary arises as an example of this construction.

Schottky groups have been generalized in multiple directions. One application of the basic ping pong argument which has seen frequent use was pioneered by Tits \cite{Pingpong}. He used what is now generally known as the Ping Pong Lemma as a tool in the proof of the famous Tits alternative. Assuming that every element in a symmetric generating set of a linear group has a unique largest eigenvalue and that they are in sufficiently generic position, he concludes that high powers of these elements generate a free group.\\
Related constructions have been used more recently to obtain nonabelian free groups acting on homogeneous spaces. Benoist used a quantitative version of transversality to define what he calls \emph{$\epsilon$-Schottky groups} acting on projective spaces \cite{BenPropre}. Using their powerful machinery of \emph{Morse actions}, Kapovich, Leeb and Porti define Schottky groups in the more general setting of semisimple Lie groups acting on symmetric spaces of noncompact type \cite{KLP}. Again, they consider a suitable notion of transversality for geodesics and show that groups generated by sufficiently high powers of axial isometries are free. In fact, they prove that such a group is Morse (or \emph{Anosov}).

A somewhat more geometric approach to generalizations of Fuchsian Schottky groups is to start out with a collection of disjoint subsets of a space $X$ and then find a set of generators pairing these subsets according to some prescribed combinatorics. Following this idea, Schottky groups were generalized to projective linear groups acting on $\mathbb{CP}^n$ (\cite{CompSchott},\cite{CompSchott2}), affine groups acting on $\mathbb{R}^3$ \cite{crookedplanes}, and the conformal Lorentzian group $\SOtt$ acting on the Einstein Universe \cite{CrookedSurf}. We introduce a new generalization which is based on the notion of cyclic orders and follows the latter approach.

The moduli space of representations $\rho:\Gamma \rightarrow G$ of a discrete group $\Gamma$ is called the $G$-\emph{character variety} of $\Gamma$. The study of character varieties and their connections to geometric structures has received a lot of attention lately. To name one example which is relevant for the paper at hand, when $\Gamma$ is the fundamental group of a closed surface and $G$ is $\PSL{2}$, there is a discrete bounded invariant called the Euler number which distinguishes connected components of the character variety \cite{BillThesis}. In this case, the component where this invariant is maximal coincides with the Teichmüller space of the surface. 

When the surface is not closed, $\Gamma$ is free and the invariant can still be defined through the use of bounded cohomology \cite{BIW}. The subset of the character variety where the Euler number is maximal is again exactly the set of representations which are holonomies of complete hyperbolic structures. The Toledo invariant generalizes this situation to representations into any Lie group $G$ of Hermitian type. The space of maximal representations in this general setting was studied by Burger, Iozzi, Wienhard, Strubel, Pozzetti (\cite{BIW},\cite{Strubel},\cite{BurgerPozzetti}). In particular, in \cite{BIW} it is shown that maximality is equivalent to the existence of an equivariant, left-continuous boundary map $\xi:\partial \bH^2 \rightarrow \Sh$, where $\Sh$ is the Shilov boundary of the symmetric space associated to $G$.

The main theorem of this paper generalizes the Schottky group construction to all maximal representations into Lie groups of Hermitian type.

In order to define the notion of Schottky group in this generality, we use \emph{partial cyclic orders} \cite{PCO}. Several notions of ordering have been successfully applied to the study of representations spaces (\cite{OrderPreserving},\cite{calegari}, \cite{GuiPositivity},\cite{GuiPositivity2}). Partial cyclic orders generalize the cyclic order on the boundary of the hyperbolic plane. The Shilov boundary of a Hermitian symmetric space admits such a partial cyclic order, defined using the generalized Maslov index introduced in  \cite{Clerc}. In the Schottky groups we introduce, pairs of points in a partially cyclically ordered space play the role of endpoints of hyperbolic geodesics. The requirement for the geodesics to be disjoint is replaced by a cyclic order condition on these pairs of points, and hyperbolic isometries are replaced by order preserving transformations. This gives a new very explicit description of maximal representations in the case of surfaces with non-empty boundary.
\begin{thmn}[\ref{SchottkyMaximal}]
Let $\Gamma = \pi_1(\Sigma)$ be the fundamental group of a compact surface with boundary, $G$ a Hermitian Lie group, and $\rho:\Gamma\to G$ a representation. Then $\rho$ is maximal if and only if it admits a Schottky presentation.
\end{thmn}

As an application of this description, we identify a domain of discontinuity for generalized Schottky subgroups of $\Sp{2n}$ acting on $\RP{2n-1}$. The key idea, suggested to us by Anna Wienhard and Olivier Guichard, is to associate certain ``halfspaces'' in $\RP{2n-1}$ to pairs of points in the Shilov boundary $\Lag{2n}$ which are used in the definition of the Schottky group. This allows us to construct a fundamental domain and analyze its orbit according to the same general strategy that one follows in the classical case of a Fuchsian Schottky group in $\PSL{2}$ acting on $\RP{1}$.

\begin{thmn}[\ref{Thm:dod}]
Let $\rho:\Gamma\to\Sp{2n}$ be a generalized Schottky group modeled on an infinite area hyperbolic surface. Then it induces a cocompact properly discontinuous action on the complement of a Cantor set of projectivized Lagrangian $n$-planes in $\RP{2n-1}$. This action admits a compact fundamental domain in $\RP{2n-1}$ bounded by finitely many smooth algebraic hypersurfaces.
\end{thmn}

In the situation of this theorem, the representation $\rho$ is Anosov and the domain of discontinuity coincides with the ones described in \cite{AnosovReps} and \cite{KLP2}.

The paper is structured as follows: In section \ref{parcyc} we introduce the notion of a partial cyclic order and some topological conditions on partially cyclically ordered sets.

Section \ref{Sec:pingpong} describes the notion of generalized Schottky group and shows the existence of an equivariant left-continuous limit curve under some topological assumptions.

In section \ref{Sec:symspaces}, we show that Shilov boundaries satisfy all the required hypotheses from the previous section and prove the main theorem as a corollary.

In the last section, we focus on the real symplectic group $\Sp{2n}$. We use our description of maximal representations as Schottky groups to construct fundamental domains bounded by algebraic hypersurfaces in $\RP{2n-1}$.

We thank William Goldman, Anna Wienhard, Olivier Guichard, and Brice Lousteau for helpful discussions and suggestions. We are grateful to Florian Stecker for his comments on a preliminary version.
\section{Partial cyclic orders} \label{parcyc}
A partial cyclic order is a relation on triples which is analogous to a partial order, but generalizing a cyclic order instead of a linear order. The definition we use was introduced in 1982 by Nov{\'a}k \cite{PCO}.
\begin{Def}
A \emph{partial cyclic order} (PCO) on a set $C$ is a relation $\cycle{}{}{}$ on triples in $C$ satisfying, for any $a,b,c,d\in C$ :
\begin{itemize}
 \item if $\cycle{a}{b}{c}$, then $\cycle{b}{c}{a}$ (\emph{cyclicity}).
 \item if $\cycle{a}{b}{c}$, then not $\cycle{c}{b}{a}$ (\emph{asymmetry}).
 \item if $\cycle{a}{b}{c}$ and $\cycle{a}{c}{d}$, then $\cycle{a}{b}{d}$ (\emph{transitivity}).
\end{itemize}
If in addition the relation satisfies:
\begin{itemize}
 \item If $a,b,c$ are distinct, then either $\cycle{a}{b}{c}$ or $\cycle{c}{b}{a}$ (\emph{totality}),
\end{itemize}
 then it is called a \emph{total cyclic order}.
\end{Def}
Let $C,D$ be partially cyclically ordered sets.
\begin{Def}
A map $f : C \rightarrow D$ is called \emph{increasing} if $\cycle{a}{b}{c}$ implies $\cycle{f(a)}{f(b)}{f(c)}$.
An automorphism of a partial cyclic order is an increasing map $f:C\rightarrow C$ with an increasing inverse. We will denote by $G$ the group of all automorphisms of $C$.
\end{Def} 
Any subset $X\subset C$ such that the restriction of the partial cyclic order is a total cyclic order on $X$ will be called a \emph{cycle}. We will also use the term cycle for (ordered) tuples $(x_1,\ldots,x_n) \in C^n$ if the cyclic order relations between the points in $C$ agree with the cyclic order given by the ordering of the tuple.
\begin{Def}
Let $a,b\in C$. The \emph{interval} between $a$ and $b$ is the set $(a,b):=\{x \in C ~|~ \cycle{a}{x}{b}\}$. The set of all intervals generates a natural topology on $C$ under which automorphisms of the partial cyclic order are homeomorphisms. We call this topology the \emph{interval topology} on $C$. We call $C$ \emph{first countable} when its interval topology is first countable.

The \emph{opposite} of an interval $I=(a,b)$ is the interval $(b,a)$, also denoted by $-I$.
\end{Def}
\begin{Ex}
The circle $S^1$ admits a total cyclic order. The relation on triples is : $\cycle{a}{b}{c}$ whenever $(a,b,c)$ are in counterclockwise order around the circle. The automorphism group of this cyclic order is the group of orientation preserving homeomorphisms of the circle.
\end{Ex}
\begin{Ex}
We can define a product cyclic order on the torus $S^1 \times S^1$. Define the relation to be $\cycle{x}{y}{z}$ whenever $\cycle{x_1}{y_1}{z_1}$ and $\cycle{x_2}{y_2}{z_2}$. This is not a total cyclic order. Some intervals in this cyclically ordered space are shown in figure \ref{fig:kthorderintervals}.
\end{Ex}
\begin{Ex}\label{Ex:linearOrder}
Every strict partial order $<$ on a set $X$ induces a partial cyclic order in the following way: define $\cycle{a}{b}{c}$ if and only if either $a<b<c$, $b<c<a$, or $c<a<b$. The cyclic permutation axiom is automatic and the two other axioms follow from the antisymmetry and transitivity axioms of a partial order.
\end{Ex}
The key topological property that we will need in the next section is a notion of completeness that we can associate to a space carrying a PCO.
\begin{Def}
A sequence $a_1,a_2,\dots \in C$ is \emph{increasing} if and only if $\cycle{a_i}{a_j}{a_k}$ whenever $i<j<k$.

Equivalently, the map $a:\mathbb{N} \rightarrow C$ defined by $a(i)=a_i$ is increasing, where the cyclic order on $\mathbb{N}$ is given by $\cycle{i}{j}{k}$ whenever $i<j<k$, $j<k<i$ or $k<i<j$ (as in Example \ref{Ex:linearOrder}).
\end{Def}
\begin{Def}
A partially cyclically ordered set $C$ is \emph{increasing-complete} if every increasing sequence converges to a unique limit in the interval topology.
\end{Def}
The following is a natural equivalence relation for increasing sequences.
\begin{Def}
Two increasing sequences $a_n$ and $b_m$ are called \emph{compatible} if they admit subsequences $a_{n_k}$ and $b_{m_l}$ making the combined sequence $a_{n_1},b_{m_1},a_{n_2},b_{m_2},\ldots$ increasing.
\end{Def}
\begin{Lem} \label{Lem:compatible_seq}
Let $C$ be an increasing-complete partially cyclically ordered set, and let $a_n$ and $b_m$ be compatible increasing sequences. Then their limits agree.
\begin{proof}
Any increasing sequence has a unique limit, and any subsequence of an increasing sequence therefore has the same unique limit.\\
The combined sequence (see the previous definition) is increasing, hence its unique limit must agree with the unique limits of both subsequences $a_{n_k}$ and $b_{m_l}$.
\end{proof}
\end{Lem}
To complete this list of definitions related to PCOs, we finish with two further restrictions on a set with a PCO which will be useful in section \ref{Sec:limitcurves}.
\begin{Def}
A partially cyclically ordered set $C$ is \emph{proper} if for any increasing quadruple $(a,b,c,d)\in C^4$, we have $\overline{(b,c)} \subset (a,d)$. Here, ``bar'' denotes the closure in the interval topology.
\end{Def}
\begin{Def}
Two points $a,b \in C$ in a partially cyclically ordered set $C$ are called \emph{comparable} if there exists a point $c\in C$ with either $\cycle abc$ or $\cycle acb$.
\end{Def}
\begin{Def}
A PCO set $C$ is \emph{full} if whenever $(a,b)$ is non-empty for some pair $a,b$, then $(b,a)$ is also non-empty. Equivalently, whenever $a,b$ are comparable then both intervals they bound are non-empty.
\end{Def}
\begin{Rem}
The motivation for the term ``full'' stems from the following construction. Assume we have a non-empty interval $(a,b)$. Then we can find a point $c\in(a,b)$ and another point $d\in(b,a)$. But then the point $d$ also lies in the interval $(c,a)$, so by fullness, there is a point inside $(a,c)$ as well. Continuing in this fashion, we can subdivide all resulting intervals further and further, and thereby construct a countably infinite subset $X\subset(a,b)$ with the following two properties: Firstly, $X$ is a cycle. Secondly, for any pair $x_1,x_2$ of distinct elements of $X$, the intersection $(x_1,x_2)\cap X$ is nonempty.
\end{Rem}
\section{Generalized Schottky groups}
\label{Sec:pingpong}
Throughout this section, $C$ denotes a partially cyclically ordered set and $G=\Aut(C)$.

\subsection{Definition of generalized Schottky group}
Let $\Sigma$ be the interior of a compact surface with boundary, of Euler characteristic $
\chi<0$. Then, the fundamental group $\pi_1(\Sigma)$ is free on $g=1-\chi$ generators. Let $\Gamma\subset \PSL{2}$ be the holonomy of a finite area hyperbolization of $\Sigma$. In this section, we construct free subgroups of $G$ using $\Gamma$ as a combinatorial model.

\begin{figure}
    \centering
    \includegraphics[width=.6\textwidth]{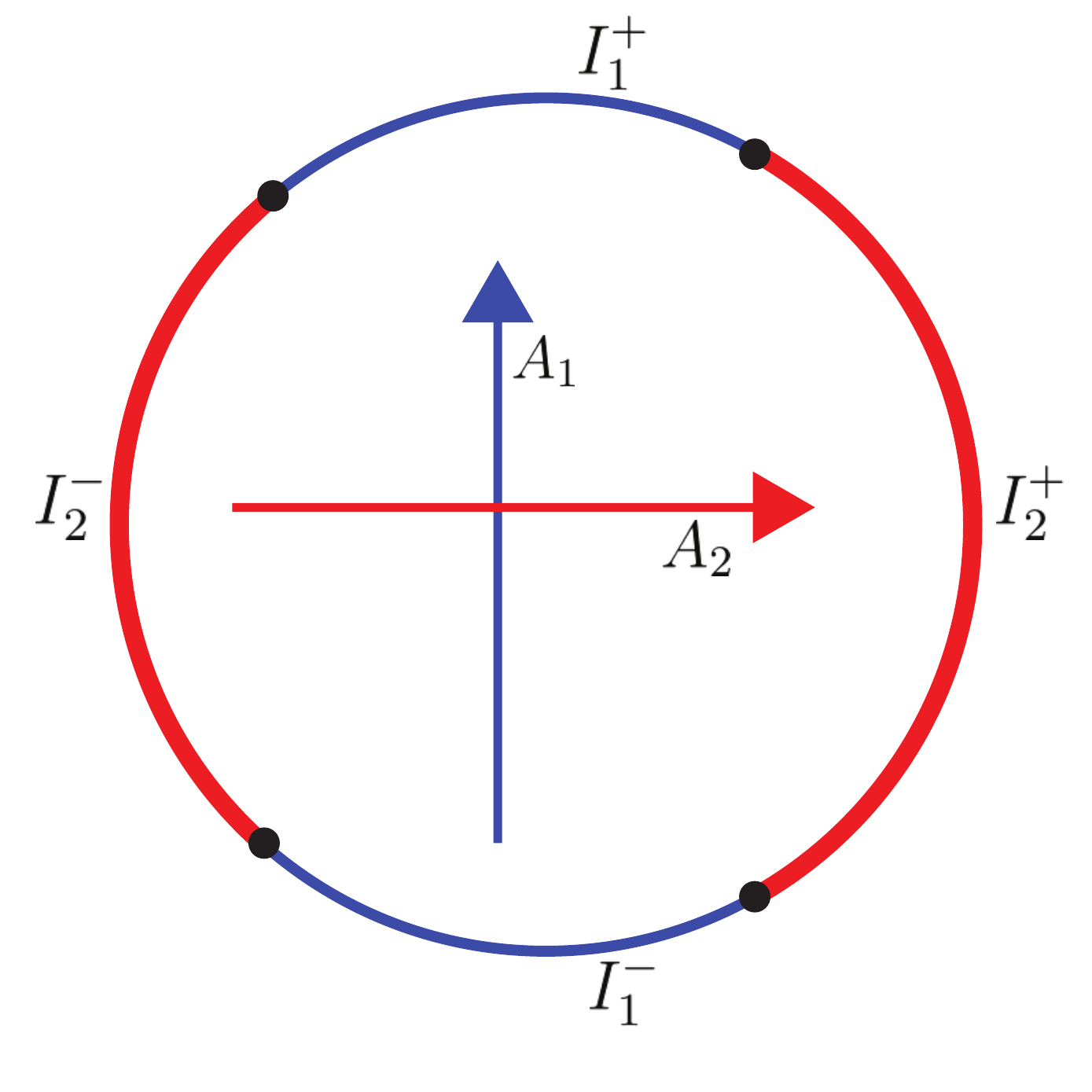}
    \caption{A combinatorial model for the once punctured torus.}
    \label{fig:Pairings}
\end{figure}

It is well known that there is a presentation for $\Gamma$ of the following form : $\Gamma$ is freely generated by $A_1,\dots,A_g \in \PSL{2}$ and there are $2g$ disjoint open intervals $I_1^+,\dots,I_g^+,I_1^-,\dots,I_g^-\subset S^1$ such that $A_j(-I_j^-)=I_j^+$. Moreover, we have that $\bigcup\limits_i \overline{I_i^+} \cup \bigcup\limits_i \overline{I_i^-} = S^1$ (figure \ref{fig:Pairings}).

The cyclic ordering on $S^1$ gives a cyclic ordering to the intervals in the definition.

We call a $k$-th order interval the image of any $I_j^+$(respectively $I_j^-$) by a reduced word of length $k-1$ not ending in $A_j^{-1}$(respectively $A_j$). There are exactly $(2g)(2g-1)^{k-1}$ $k$-th order intervals. There is a natural bijection between words of length $k$ and $k$-th order intervals. We use this bijection to index $k$-th order intervals : $I_W$ is the interval corresponding to the word $W$. For any fixed $k$, the $k$-th order intervals are all pairwise disjoint, and so they are cyclically ordered. This induces a cyclic ordering on words of length $k$ in $\Gamma$. The union of all closures of $k$-th order intervals is all of $S^1$. 

The following easy Lemma, which is a reformulation of transivity, motivates our definition of generalized Schottky groups in $G$.
\begin{Lem}
Let $(a,b,c)\in C^3$ be a cycle. Then we have $(b,c) \subset (b,a)$. In particular, the intervals $(a,b)$ and $(b,c)$ are disjoint.
\begin{proof}
Let $x\in(b,c)$, so we have $\cycle bxc$. By transitivity, together with $\cycle bca$, this implies $\cycle bxa$.
\end{proof}
\end{Lem}
We now define generalized Schottky subgroups of $G$ by asking for a setup of intervals similar to the $\PSL{2}$ case and requiring generators to pair them the same way.
\begin{Def}
Let $\xi_0$ be an increasing map from the set of endpoints of the intervals $I_1^+,\dots,I_g^+$, $I_1^-,\dots,I_g^-$ into $C$. For $I_i^\pm=(a_i^\pm,b_i^\pm)$, define the corresponding interval $J_i^\pm=(\xi_0(a_i^\pm),\xi_0(b_i^\pm))\subset C$. Next, assume there exist $h_1,\dots,h_g\in G$ which pair the endpoints of $J_i^\pm$ in the same way that the $g_i$ pair the $I^\pm_i$, so that $h_i(-J_i^-)=J_i^+$. We call the image of the induced map $\Gamma \to G$ sending $A_i$ to $h_i$ a \emph{generalized Schottky group}, and the intervals $J_i^\pm$ used to define it a set of \emph{Schottky intervals} for this group.
\end{Def}
\begin{Rem}\
\begin{enumerate}
    \item A generalized Schottky group will in general have many possible choices of a set of Schottky intervals. We will only use this term when a specific choice of both generators and intervals corresponding to these generators is fixed.
    \item Since the cyclic ordering is a property of $\RP{1}$ which is not shared by $\mathbb{CP}^1$, the Schottky groups defined here do \emph{not} generalize the more well known $\mathbb{CP}^1$ Kleinian case.
    \item The requirement that the combinatorial model be a finite-volume hyperbolization is artificial. It is helpful in order to avoid having to separate our analysis into several cases. We could use a model where the intervals have disjoint closures and the construction would work in the same way. In section \ref{sec:FundamentalDomains} we will use presentations of this type to describe domains of discontinuity in $\RP{2n-1}$. 
    \item Our use of the term "Schottky" differs slightly from most references in that we allow for the closures of the ping pong subsets to intersect. This is sometimes called "Schottky-type".
\end{enumerate}
\end{Rem}

\begin{figure}
    \centering
    \includegraphics[width=.6\textwidth]{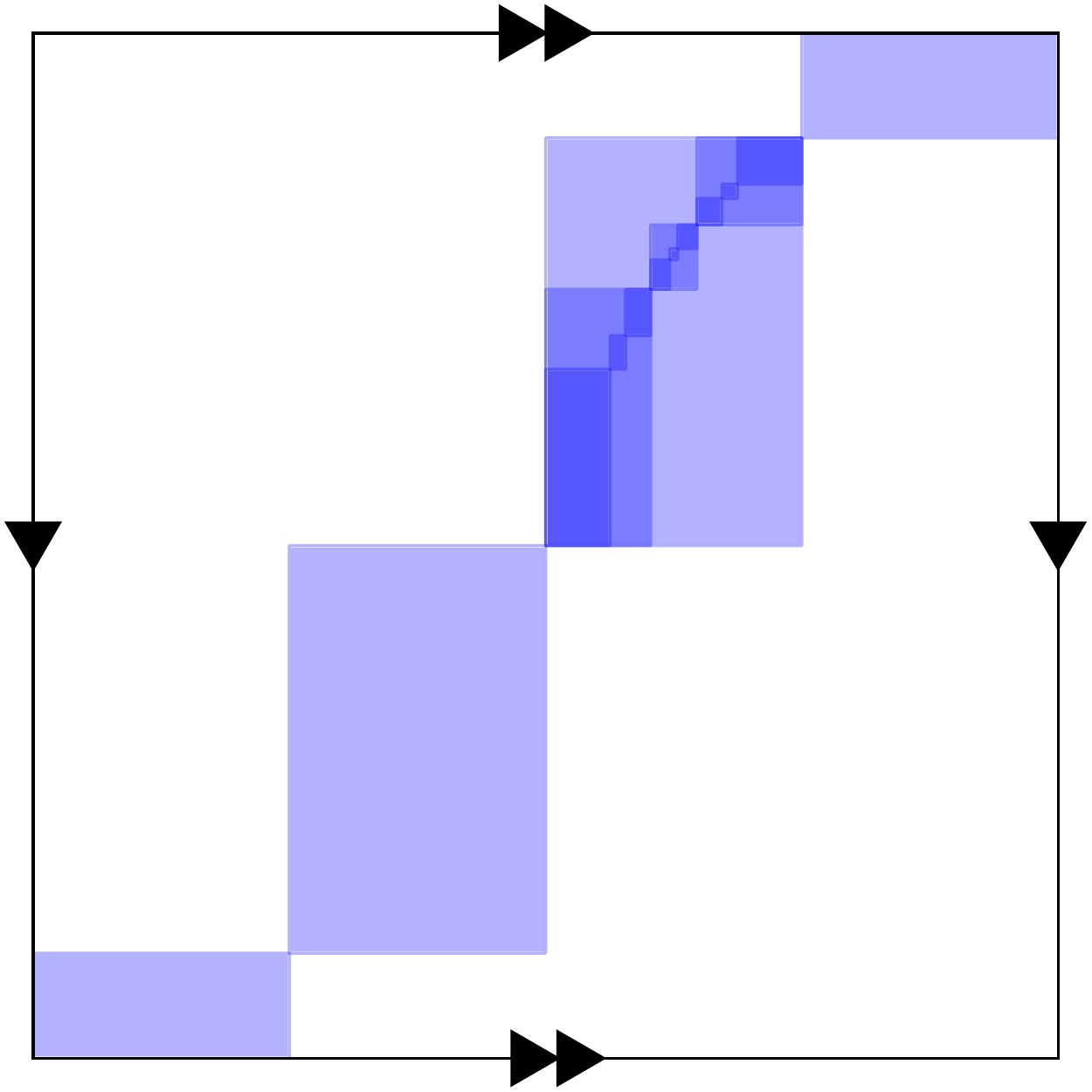}
    \caption{Some first, second and third order intervals for a generalized Schottky group acting on $S^1 \times S^1$.}
    \label{fig:kthorderintervals}
\end{figure}

With this setup, we can define $k$-th order intervals in $C$ in the same way as above but starting with the intervals $J_i^\pm$ and their images under words in the $h_i$ (see figure \ref{fig:kthorderintervals}). As above, denote by $J_W$ the interval corresponding to $W$. Note that since $\xi_0$ is increasing, the $k$-th order intervals in $C$ are also cyclically ordered, where the ordering is the same as the ordering of the corresponding intervals in $S^1$.
\begin{Prop}
The group generated by $h_1\dots h_g$ is free on those generators.
\begin{proof}
Define $J_i = J_i^+ \cup J_i^-$. Note that $J_i\cap J_j = \emptyset$ whenever $i\neq j$. Moreover, for any $n\neq 0$, $h_i^n (J_j) \subset J_i$ and so the proposition follows from the ping pong lemma.
\end{proof}
\end{Prop}

The endpoints of $k$-th order intervals in $C$ satisfy the same cyclic order relations as the corresponding endpoints in $S^1$, and we can extend $\xi_0$ to an increasing equivariant map defined on the countable dense set of all endpoints of Schottky intervals in $S^1$. We denote this set by $S^1_\Gamma$.
\subsection{Limit curves}\label{Sec:limitcurves}
\begin{Lem} \label{Lem:closures_pos}
Let $C$ be a partially cyclically ordered set which is full and proper, and $(x_1,\ldots,x_6) \in C^6$ a cycle. Let $I_1 = (x_1,x_2), I_2=(x_3,x_4),I_3=(x_5,x_6)$ and $a_i \in \overline{I_i}$ be arbitrary points in the closures of the intervals. Then we have $\cycle {a_1}{a_2}{a_3}$.
\begin{proof}
Since $C$ is full, we can choose auxiliary points $y_i$ such that the $12$-tuple $(y_1,x_1,x_2,y_2,y_3,x_3,x_4,y_4,y_5,x_5,x_6,y_6)$ is a cycle. This allows us to conclude that $\overline{(x_i,x_{i+1})} \subset (y_{i},y_{i+1})$ for odd $i$ as $C$ is proper. Since $(y_1,\ldots,y_6)$ is a cycle, transitivity implies the lemma.
\end{proof}
\end{Lem}
\begin{Lem}
\label{Lem:nested_sequences}
Let $P_n \rightarrow P$ be an increasing sequence in a proper, increasing complete, PCO set $C$. Assume $Q_n$ is another sequence with $Q_n\in \overline{(P_n,P_{n+1})}$ for all $n$. Then $Q_n$ converges to $P$ and is $3$-increasing in the following sense: whenever $i+2<j<k-2$, we have $\cycle{Q_i}{Q_j}{Q_k}$.
\begin{proof}
For every $n\geq 2$, $Q_n\in (P_{n-1},P_{n+2})$ by properness, which already implies that $Q_n$ is $3$-increasing. Now, consider the following sequence: \[P_1,Q_2,P_4,Q_5,\dots,P_{3n+1},Q_{3n+2},\dots\] 
It is increasing, and admits a subsequence which is also a subsequence of $P_n$. Since increasing sequences have unique limits, this sequence must converge to $P$. The increasing subsequence $Q_{3n+2}$ therefore converges to $P$. Using the same argument, we see that $Q_{3n+1}$ and $Q_{3n}$ also converge to $P$, so in fact the sequence $Q_n$ converges to $P$.
\end{proof}
\end{Lem}
We now come to the main theorem of this section, which explains how to construct a boundary map for generalized Schottky groups, under some topological assumptions.\\
The boundary map we construct will be left-continuous as a map from $S^1$ to some first countable topological space $C$. To avoid confusion, let us fix the definition here: In a small neighborhood $U$ of a point $a\in S^1$, the cyclic order induces a linear order. A sequence $a_n \in U$ converges to $a$ from the left if $a_n<a$ and $a_n\xrightarrow{n\to\infty}a$. The function $f$ is left-continuous at $a$ if $f(a_n)\to f(a)$ for all sequences $a_n$ converging to $a$ from the left.\\
It is worth noting that to check for left-continuity at a point $a$, it is in fact sufficient to check the convergence of $f(a_n)$ for increasing sequences $a_n$ converging to $a$ from the left. The reason is the following: Assume $a_n$ is a sequence converging to $a$ from the left such that $f(a_n)$ does not converge to $f(a)$. Then it has a subsequence such that $f(a_{n_k})$ stays bounded away from $f(a)$. But since $a_{n_k} \to a$ from the left, we can pick a further subsequence which is increasing and still a counterexample to left-continuity.
\begin{Thm}    \label{Thm:bdry_map_gen_Schottky}
Let $\rho: \Gamma \to G$ be the map defining a generalized Schottky group. Assume that $C$ is first countable, increasing-complete, full and proper. Then there is a left-continuous, equivariant, increasing boundary map $\xi: S^1 \to C$.
\begin{proof}
We construct the map $\xi$ as follows: Recall that $S^1_\Gamma \subset S^1$ denotes the domain of $\xi_0$ and is a dense subset. For $x\in S^1$, pick any increasing sequence $x_n \in S^1_\Gamma$ converging to $x$ and set
\[ \xi(x) = \lim_{n\to\infty} \xi_0(x_n). \]
First of all, let us show that this value is well-defined. Since $x_n$ is an increasing sequence, the increasing map $\xi_0$ maps it to an increasing sequence in $C$ which therefore has a unique limit. Furthermore, this limit does not depend on the choice of $x_n$: Let $y_m$ be another increasing sequence converging to $x$. Then the two sequences $\xi_0(x_n)$ and $\xi_0(y_m)$ are compatible, so they have the same limit by Lemma \ref{Lem:compatible_seq}.

We then verify that $\xi$ is equivariant. Let $x\in S^1$, $\gamma\in\Gamma$, and $x_n\to x$ an increasing sequence, so we have $\xi(x) = \lim\limits_{n\to\infty} \xi_0(x_n)$. Then $\gamma(x_n)$ is an increasing sequence converging to $\gamma(x)$, so by continuity of $\rho(\gamma)$ and equivariance of $\xi_0$, we have the following equalities:
\[ \rho(\gamma)(\xi(x)) = \lim_{n\to\infty} \rho(\gamma)(\xi_0(x_n)) = \lim_{n\to\infty} \xi_0(\gamma(x_n)) = \xi(\gamma(x)). \]

Next, we show that it is left-continuous. Assume $x_n \in S^1$ is a sequence converging to $x$ from the left. As explained above, without loss of generality we can take $x_n$ to be an increasing sequence. We pick points $y_n \in S^1_\Gamma$ such that $y_n \in (x_{n-1},x_n)$. Then $y_n$ is increasing and we have $x_n\in (y_n,y_{n+1})$. Furthermore, $y_n$ also converges to $x$, hence 
\begin{equation}    \label{bdry_map_eq1}
 \xi(x) = \lim \xi_0(y_n).
\end{equation}
Now, for each $n$, let $\{a_k(n)\}_{k\in\bN}\subset S^1_\Gamma$ be an increasing sequence converging to $x_n$, so
\begin{equation}    \label{bdry_map_eq2}
 \xi(x_n) = \lim_{k\to\infty} \xi_0(a_k(n)).
\end{equation} 
Then $a_k(n) \in (y_n,y_{n+1})$ for large $k$, so 
\begin{equation}    \label{bdry_map_eq3}
 \lim\limits_{k\to\infty} \xi_0(a_k(n)) \in \overline{(\xi_0(y_n),\xi_0(y_{n+1}))}
\end{equation} 
because $\xi_0$ is increasing. Now Lemma \ref{Lem:nested_sequences} applies and, combined with \eqref{bdry_map_eq1}, \eqref{bdry_map_eq2} and \eqref{bdry_map_eq3}, tells us that $\xi(x_n)$ converges to $\xi(x)$. 

The final property we need to check is that $\xi$ is increasing. Assume that we have $\cycle xyz$ for points $x,y,z \in S^1$. By density of $S^1_\Gamma$, we can find a cycle $(a_1,a_2,b_1,b_2,c_1,c_2) \in (S^1_\Gamma)^6$ such that $x\in (a_1,a_2), \ y\in (b_1,b_2), \ z\in (c_1,c_2)$. As in the proof of left-continuity, this implies that $\xi(x) \in \overline{(\xi_0(a_1),\xi_0(a_2))}$, and similar for the other two points. Using Lemma \ref{Lem:closures_pos}, we conclude $\cycle{\xi(x)}{\xi(y)}{\xi(z)}$.
\end{proof}
\end{Thm}
The very general construction described in this section applies to many examples. For instance, various notions of positivity in homogeneous spaces give rise to partial cyclic orders. More specifically, the Shilov boundary of Hermitian symmetric spaces admits a PCO satisfying all the above properties, and the rest of the paper is dedicated to this example. It is also possible, using techniques similar to Fock-Goncharov total positivity \cite{FocPositivity}, to construct a PCO on spaces of complete oriented flags. In this last case, for flags in $\mathbb{R}^3$, Schottky groups can be used to describe convex projective structures on surfaces with boundary \cite{NextPaper}.
\section{Hermitian symmetric spaces of tube type}\label{Sec:symspaces}
In this section, we show that the Shilov boundary of a Hermitian symmetric space of tube type $X$ admits a partial cyclic order invariant under the holomorphic isometry group of $X$. Moreover, we prove that Shilov boundaries satisfy the topological assumptions from theorem \ref{Thm:bdry_map_gen_Schottky}, so we have a boundary map for every generalized Schottky subgroup. Then, using the machinery of section \ref{Sec:pingpong}, we show that Schottky subgroups in this case correspond to maximal representations.

\subsection{The partial cyclic order on the Shilov boundary}
Consider a tube type domain $X=V+i\Omega$, where $V$ is a real vector space and $\Omega\subset V$ is a symmetric cone. Let $G$ be the group of holomorphic isometries of $X$.
\begin{Def} The \emph{Shilov boundary} $\Sh$ of $X$ is the unique closed $G$-orbit in the visual boundary of $X$. It can be identified with a subset of the topological boundary of the bounded domain realization of $X$. The action of $G$ extends smoothly to $\Sh$.
\end{Def}
The vector space $V$ admits a Jordan algebra structure associated to the symmetric cone $\Omega$ \cite{farautkoranyi}.
\begin{Def} The partial order $<_\Omega$ on $V$ is given by $x<_\Omega y$ if and only if $y-x\in\Omega$. Equivalently, $y-x = z^2$ for some $z\in V$.\end{Def}
\begin{Prop}[\cite{farautkoranyi}, Proposition X.2.3]
The \emph{Cayley Transform} $p:V \rightarrow V_\mathbb{C}$ defined by
\[p(v) = (v-i e)(v+i e)^{-1},\]
where $e$ is the identity of the Jordan algebra, maps the vector space $V$ into the Shilov boundary $\Sh$ and $\overline{p(V)}=\Sh$.
\end{Prop}
Using the notion of transversality, we can make the statement more precise.
\begin{Def}
Two points $x,y\in \Sh$ are called \emph{transverse} if the pair $(x,y)\in \Sh\times \Sh$ belongs to the unique open $G$-orbit for the diagonal action.
\end{Def}
The image of the Cayley transform is exactly the set of points $x\in \Sh$ which are transverse to a fixed point which we denote by $\infty$. \cite[Section 6.6.1]{AnnaThesis}

The next object we need to define is the generalized Maslov index. This index is a function on ordered triples of points in $\Sh$, invariant under $G$. It will be used in order to define a partial cyclic order on $\Sh$, extending the partial cyclic order induced by $<_\Omega$ on $p(V)\subset \Sh$.

The generalized Maslov index is defined in \cite{Clerc} using the notion of $\Gamma$-radial convergence. For our purposes we will use the following equivalent definition, given in the same paper.
\begin{Def}
Let $x,y,z\in \Sh$. Applying an element of $G$, we may assume $x,y,z\in p(V)$. Let $v_x,v_y,v_z\in V$ be the vectors which map respectively to $x,y,z$ under the Cayley transform $p$. Then, the \emph{generalized Maslov index} of $x,y,z$ is the integer
\[\M(x,y,z) := \mathsf{k}(v_y-v_x) + \mathsf{k}(v_z-v_y) + \mathsf{k}(v_x-v_z)\]
where $\mathsf{k}(v)$ is the difference between the number of positive eigenvalues of $v$ and the number of negative eigenvalues of $v$ in its spectral decomposition.

When $x,y$ are transverse to $z$, equivalently, we can map $z$ to $\infty$ using an element of $G$ and define
\[\M(x,y,\infty) = \mathsf{k}(v_y-v_x)\] 
\end{Def}
\begin{Prop}
The Maslov index enjoys the following properties:
\begin{itemize}
\item $G$-invariance : $\M(gx,gy,gz) = \M(x,y,z)$.
\item Skew-symmetricity : $\M(x_1,x_2,x_3) = \sgn(\sigma)\M(x_{\sigma(1)},x_{\sigma(2)},x_{\sigma(3)}).$
\item Cocycle identity : $\M(y,z,w) - \M(x,z,w) + \M(x,y,w) - \M(x,y,z) = 0$.
\item Boundedness : $|\M(x,y,z)| \leq \mathrm{rk}(X)$ 
\end{itemize}
\end{Prop}
These properties allow us to define a partial cyclic order on the Shilov boundary.
\begin{Prop}
The relation $\cycle{x}{y}{z}$ if and only if $\M(x,y,z)=\mathrm{rk}(X)$ defines a $G$-invariant partial cyclic order on $\Sh$.
\begin{proof}
Since $M$ is skew-symmetric, the relation automatically satisfies the first two axioms of a partial cyclic order. To prove the third axiom, assume $\M(x,y,z)=\M(x,z,w)=\mathrm{rk}(X)$. By the cocycle identity,
\[\M(y,z,w) - \M(x,z,w) + \M(x,y,w) - \M(x,y,z) = 0\]
and so
\[\M(y,z,w) + \M(x,y,w) = 2rk(X)\]
which is only possible if $\M(Q,R,S)=\M(P,Q,S)=\mathrm{rk}(X)$.
\end{proof}
\end{Prop}
The partial cyclic order $\cycle{}{}{}$ is closely related with the causal structure on $\Sh$ introduced by Kaneyuki \cite{KanShilov}. Namely, whenever $\cycle{x}{y}{z}$, there is a future-oriented timelike curve going through $x,y,z$ in that order. Informally, $y$ is in the intersection of the future of $x$ and the past of $z$. The following two lemmas describe some immediate properties of cyclically ordered triples.
\begin{Lem}[\cite{AnnaThesis}, Lemma 5.5.4] \label{Lem:PosTriplesTransverse}
Let $x,y,z \in \Sh$ with $\cycle xyz$. Then $x,y,z$ are pairwise transverse.
\end{Lem}
\begin{Lem} \label{Lem:Intervals_V}
Assume $x,y\in V$. Then, $\cycle{x}{y}{\infty}$ if and only if $x<_\Omega y$.
\begin{proof}
The cone $\Omega$ coincides with the region where $\mathsf{k}(v)=\mathrm{rk}(X)$.
\end{proof}
\end{Lem}
\begin{Rem}
The interval topology on $\Sh$ is the same as the usual manifold topology.
\end{Rem}
\begin{Prop} \label{Thm:Shilov_complete}
The PCO defined by $\cycle{}{}{}$ on $\Sh$ is increasing-complete, full and proper.
\begin{proof}
We first show that it is increasing-complete. Let $x_1,x_2,\dots$ be an increasing sequence in $\Sh$. Let $g\in G$ be such that $gx_2 = \infty$. Then, since we have $\cycle{x_k}{x_{k+1}}{x_2}$ for all $k\geq 3$, the sequence $gx_3,gx_4,\dots$ is an increasing sequence transverse to $\infty$. Hence, there exist $v_3,v_4,\ldots \in V$ with $p(v_k)=gx_k$.

This new sequence is increasing with respect to $<_\Omega$. Moreover, it is bounded since we have $\cycle{gx_k}{gx_1}{gx_2}$ for all $k>2$, so $v_k<_\Omega v_1$ where $p(v_1)=gx_1$. The tail of the sequence is contained in $\overline{(v_3,v_1)}$ which is compact, so it has an accumulation point. If $w,w'$ are two accumulations points of the sequence, let $w_k$, $w'_k$ be subsequences converging respectively to each of them. Passing to subsequences if necessary, we can arrange so that $w_k<_\Omega w'_k$ for all $k$, and so $w'_k-w_k\in \Omega$. This implies $w'-w\in \overline{\Omega}$, and by the same argument we can also show $w-w'\in \overline{\Omega}$. Since $\Omega$ is a proper convex cone (in the sense of \cite{farautkoranyi}), its closure does not contain any opposite pairs, so $w=w'$.

Now we turn to fullness of the PCO. Whenever an interval $(x,y)$ is nonempty, its endpoints have to be transverse by Lemma \ref{Lem:PosTriplesTransverse}. We can therefore apply an element of $G$ to map $x$ to $\infty$ and $y$ inside $p(V)$. Then Lemma \ref{Lem:Intervals_V} shows that the interval $(y,x)$ is also nonempty.

Finally, we show that the PCO is proper. Let $(x_1,x_2,x_3,x_4) \in \Sh^4$ be a cycle. Using an element of $G$, we can assume that $x_4$ is $\infty$, so that $x_1,x_2,x_3 \in p(V)$. Let $v_i \in V$ be the vector such that $p(v_i) = x_i$ for $i=1,2,3$. Now the cyclic relations $\cycle{x_1}{x_2}{\infty}$ and $\cycle{x_2}{x_3}{\infty}$ imply that both $v_2-v_1$ and $v_3-v_2$ lie in the cone $\Omega$. The interval $(x_2,x_3)$ is therefore given by $p\left( (v_2 + \Omega) \cap (v_3 - \Omega)\right)$. This implies the claim since $(v_2 + \Omega) \cap (v_3 - \Omega)$ is a relatively compact set in $V$ whose closure is contained in $v_1 + \Omega$, which is mapped onto $(x_1,\infty)$ by $p$.
\end{proof}
\end{Prop}
\subsection{Maximal representations}
In the previous section we defined a PCO on the Shilov boundary $\Sh$ of a Hermitian symmetric space of tube type on which the group of isometries $G$ acts by order-preserving diffeomorphisms. We recall that this action is transitive on transverse pairs. The Schottky construction described in section \ref{Sec:pingpong} therefore gives maps $\rho:\Gamma\rightarrow G$ where $\Gamma$ is the fundamental group of a surface with boundary. Maximal representations, defined and studied in \cite{BIW}, are a class of geometrically interesting representations and we will show in this section that they correspond to Schottky subgroups.

The key is the following characterization from \cite{BIW} (Theorem 8) :
\begin{Thm}
Let $h:\Gamma\rightarrow \PSL{2}$ be a complete finite area hyperbolization of the interior of $\Sigma$ and $\rho:\Gamma\rightarrow G$ a representation into a group of Hermitian type. Then $\rho$ is maximal if and only if there exists a left continuous, equivariant, increasing map
\[\xi : S^1 \rightarrow \Sh\]
where $\Sh$ is the Shilov boundary of the bounded symmetric domain associated to $G$.
\end{Thm}
Using this characterization and our earlier construction of a boundary map for generalized Schottky representations, we see that the two notions agree:
\begin{Thm}\label{SchottkyMaximal}
The representation $\rho:\Gamma\rightarrow G$ is maximal if and only if it admits a Schottky presentation.
\begin{proof}
Assume $\rho$ is Schottky. Proposition \ref{Thm:Shilov_complete} states that all the prerequisites of Theorem \ref{Thm:bdry_map_gen_Schottky} are fulfilled. Therefore, there exists a boundary map $\xi$ satisfying the conditions of the characterization above, so $\rho$ is maximal.

Conversely, if $\rho$ is maximal, then we have such a map $\xi$. Choosing a Schottky presentation for the hyperbolisation $h$, we get a Schottky presentation for $\rho$ by using the intervals $(\xi(a),\xi(b))$ where $(a,b)$ is some Schottky interval in the presentation for $h$. Equivariance and positivity of $\xi$ ensure that these intervals fit our definition of generalized Schottky groups.
\end{proof}
\end{Thm}
Theorem \ref{SchottkyMaximal}, as stated, assumes that $G$ is of tube type. However, this assumption is not necessary. This is because of the following observations. Let $X$ be a Hermitian symmetric space, and $\Sh$ its Shilov boundary. Then, in the same way as for tube type, the generalized Maslov index defines a partial cyclic order on $\Sh$. Let $x,y\in \Sh$ be transverse. Then, $x,y$ are contained in the Shilov boundary of a unique maximal tube type subdomain of $X$ \cite[Lemma 4.4.2]{AnnaThesis}. Moreover, this is also true of any increasing triple in $\Sh$ \cite[Proposition 5.1.4]{AnnaThesis}. This means that any increasing subset of $\Sh$ is contained in the Shilov boundary of a tube type subdomain, and so the proofs of this section generalize to arbitrary Hermitian symmetric spaces.
\section{Schottky groups in $\Sp{2n}$}

In this section, we consider the symplectic group $\Sp{2n}$, acting on $\bR^{2n}$ equipped with a symplectic form $\omega$, and describe the construction of Schottky groups in detail.
\subsection{The Maslov index in $\Sp{2n}$}
\begin{Def}
Let $P,Q$ be transverse Lagrangians in $\bR^{2n}$. We associate to them an antisymplectic involution $\Inv{P}{Q}$ defined using the splitting $\bR^{2n}=P\oplus Q$:
\begin{align*} \Inv{P}{Q}: P\oplus Q & \to P\oplus Q \\
 (v,w) & \mapsto (-v,w)
 \end{align*}
We call this antisymplectic involution the \emph{reflection} in the pair $P,Q$. This generalizes the projective reflection in $\RP{1}$.

We will sometimes abuse notation and use $\Inv{P}{Q}$ to denote the induced transformation on Grassmannians.
\end{Def}
Using this involution, we associate a symmetric bilinear form to the pair $P,Q$ :
\begin{Def}
\[\B{P}{Q}(v,w) := \omega(v, \Inv{P}{Q}(w))\]
\end{Def}
This bilinear form is nondegenerate and has signature $(n,n)$.
\begin{Def} \label{Def:MaslovIndex}
Let $P,Q, R$ be pairwise transverse Lagrangians in $\bR^{2n}$. The \emph{Maslov index} of the triple $(P,Q,R)$ is the index of the restriction of $\B{P}{R}$ to $Q$. We denote it by $\M(P,Q,R)$.
\end{Def}

Since $\Lag{2n}$ is the Shilov boundary for the bounded domain realization of the symmetric space of $\Sp{2n}$, it is an example of the general construction in section \ref{Sec:symspaces}. In fact, the Maslov index we just defined agrees with the more general version that we introduced before. Hence, the relation defined by $\cycle{P}{Q}{R}$ whenever $\M(P,Q,R)=n$ is a partial cyclic order on $\Lag{2n}$, enabling us to apply the constructions and results from section $\ref{Sec:pingpong}$.\\
We also remark that the definition makes sense for any isotropic subspace $Q$, not only the maximal isotropic ones. 

The following property of the Maslov index is well-known.

\begin{Prop}
    The Maslov index classifies orbits of triples of pairwise transverse Lagrangians, i.e. the map 
    \[ (P,Q,R) \mapsto \M(P,Q,R) \] induces a bijection from orbits of pairwise transverse Lagrangians under $\Sp{2n}$ to the set $\{ -n,-n+2,\ldots,n\}$.
\end{Prop}
The Maslov index and the reflection in a pair of Lagrangians are related in the following way:
\begin{Prop}
\[\M(P,\Inv{P}{Q}(V),Q)=-\M(P,V,Q).\]
\begin{proof}
\[\B{P}{Q}\big(\Inv{P}{Q}(u),\Inv{P}{Q}(v)\big)=\omega(\Inv{P}{Q}(u),v)=-\omega(u,\Inv{P}{Q}(v))=-\B{P}{Q}(u,v). \qedhere  \]
\end{proof}
\end{Prop}
 The proposition above means that reflections reverse the partial cyclic order.
\subsection{Fundamental domains}
\label{sec:FundamentalDomains}
In the special case of $\Sp{2n}$, the Schottky groups we obtain admit nice fundamental domains for their action on $\RP{2n-1}$. The domain of discontinuity which is the orbit of this fundamental domain is in general hard to describe, but it simplifies in some cases.\\
We will proceed as follows: First, we associate a ``halfspace'' in $\RP{2n-1}$ to each interval in $\Lag{2n}$ and explain how to construct the fundamental domain. Then we cover some preliminaries which will allow us to explicitly identify the domain of discontinuity for generalized Schottky groups modeled on an infinite volume hyperbolization. More specifically, we explain how to identify an interval with the symmetric space associated with $\GL{n}$ and how to use a contraction property from \cite{Contraction} for maps sending one interval into another.

%
\subsubsection{Positive halfspaces and fundamental domains}\label{sec:pos_halfspaces}
\begin{Def}
Let $P,Q$ be an ordered pair of transverse Lagrangians. We define the $\emph{positive halfspace}$ $\PH{P}{Q}$ as the subset
\[\PH{P}{Q}:=\{ \ell \in \RP{2n-1} ~|~ \B{P}{Q}|_{\ell\times\ell} > 0 \}.\]
It is the set of positive lines for the form $\B{P}{Q}$.
\end{Def}
The positive halfspace $\PH{P}{Q}$ is bounded by the conic defined by $\B{P}{Q}=0$. This type of bounding hypersurface was introduced by Guichard and Wienhard in order to describe Anosov representations of closed surfaces into $\Sp{2n}$. They are also the boundaries of $\mathbb{R}$-tubes defined in \cite{BurgerPozzetti}. A symplectic linear transformation $T\in\Sp{2n}$ acts on positive halfspaces in the following way : $T\PH{P}{Q} = \PH{TP}{TQ}$.

\begin{figure}
    \centering
    \includegraphics[width=\textwidth]{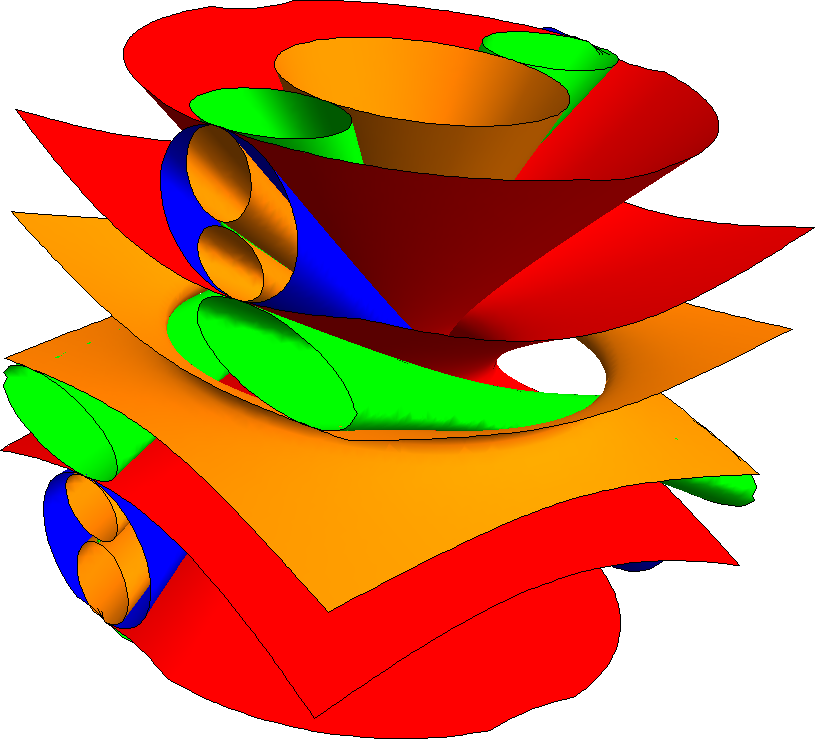}
    \caption{The first two generations of positive halfspaces for a two-generator Schottky group in $\Sp{4}$.}
    \label{fig:sp4halfspaces}
\end{figure}

\begin{Prop}
Let $P,Q$ be an ordered pair of Lagrangians. Then,
\[\PH{Q}{P}=\overline{\PH{P}{Q}}^C = \Inv{P}{Q}(\PH{P}{Q})\]
\begin{proof}
For the first equality,
\[\B{Q}{P}(v,w)=\omega(v,\Inv{Q}{P}(w))=\omega(v,-\Inv{P}{Q}(w))=-\B{P}{Q}(v,w).\]
For the second equality, notice that $\B{P}{Q}(\Inv{P}{Q}(v),\Inv{P}{Q}(w)) = -\B{P}{Q}(v,w)$.
\end{proof}
\end{Prop}
\begin{Prop}
\label{Prop:proj}
A positive halfspace is the projectivisation of an interval, that is,
\[\PH{P}{Q}=\bigcup_{L\in (P,Q)} \mathbb{P}(L)\]
\begin{proof}
If $\ell \subset L$ for some $L\in(P,Q)$, then
\[\B{P}{Q}|_{\ell\times\ell} > 0\]
and so $\ell \in \PH{P}{Q}$.

Conversely, if $\ell \in \PH{P}{Q}$, then we wish to find a Lagrangian $L\supset \ell$ with $\M(P,L,Q)=n$. Consider the subspace $V=\langle\ell,\Inv{P}{Q}(\ell)\rangle$. The form $\B{P}{Q}$ has signature $(1,1)$ on that subspace, and so its orthogonal has signature $(n-1,n-1)$. Moreover, the form $\omega$ is nondegenerate on $V$ so $V^{\perp_\omega}$ is a symplectic subspace. Notice that
\[V^{\perp_\BB}=\langle\ell,\Inv{P}{Q}(\ell)\rangle^{\perp_\BB} = \ell^{\perp_\BB} \cap (\Inv{P}{Q}(\ell))^{\perp_\BB}= \ell^{\perp_\BB} \cap \ell^{\perp_\omega}=V^{\perp_\omega}.\]
So we can pick a positive definite Lagrangian $L'\subset V^\perp$, which will be orthogonal to $\ell$ for both $\omega$ and $\B{P}{Q}$, so $L=\langle L',\ell \rangle$ is a positive definite Lagrangian containing $\ell$.
\end{proof}
\end{Prop}
\begin{Lem}
\label{Lem:Disjointness}
If $(P,Q,R,S)$ is a cycle in $\Lag{2n}$ and $V\in\Lag{2n}$ such that $\M(P,V,Q)=n$, then $\M(R,V,S)=-n$.
\begin{proof}
Using the cocycle relation,
\[\M(V,Q,R)-\M(P,Q,R)+\M(P,V,R)-\M(P,V,Q)=0\]
so \[\M(V,Q,R)+\M(P,V,R)=2n\]
which implies that $\M(V,Q,R)=\M(P,V,R)=n$.

Similarly,
\[\M(V,R,S)-\M(P,R,S)+\M(P,V,S)-\M(P,V,R)=0\]
so
\[\M(V,R,S)+\M(P,V,S)=2n\]
which means that $\M(V,R,S)=n$ and so $\M(R,V,S)=-n$.
\end{proof}
\end{Lem}
Now we can prove the disjointness criterion for positive halfspaces.
\begin{Prop}
\label{Prop:Disjointness}
If $(P,Q,R,S)$ is a cycle in $\Lag{2n}$, then $\PH{P}{Q}$ is disjoint from $\PH{R}{S}$.
\begin{proof}
Let $\ell\in \PH{P}{Q}$. By Proposition \ref{Prop:proj}, $\ell \subset L$ for some Lagrangian $L$ with $\M(P,L,Q)=n$. By Lemma \ref{Lem:Disjointness}, $\M(R,L,S)=-n$ which means that $\B{R}{S}|_\ell<0$ and so $\ell \notin \PH{R}{S}$.
\end{proof}
\end{Prop}
\begin{figure}
    \centering
    \includegraphics[width=.4\textwidth]{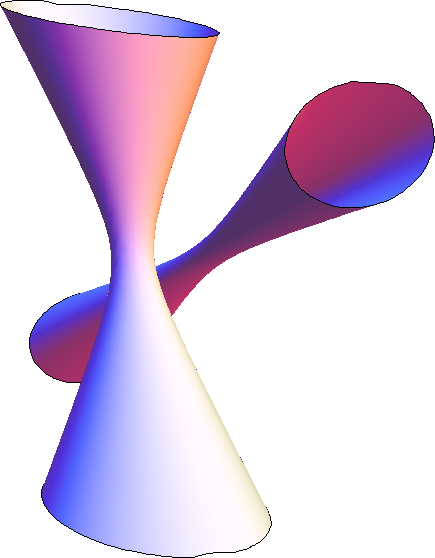}
    \caption{A pair of disjoint positive halfspaces in $\RP{3}$}
    \label{fig:halfspaces}
\end{figure}
For any generalized Schottky group, we can use this previous proposition to construct a fundamental domain. If the defining intervals for the Schottky group are $(a_1^\pm,b_1^\pm),\dots,(a_g^\pm,b_g^\pm)\subset \Lag{2n}$, let
\[D = \bigcap_{j=1}^g \left(\PH{a_j^+}{b_j^+}\cup\PH{a_j^-}{b_j^-}\right)^C.\]
That is, $D$ is the subset of $\RP{2n-1}$ which is the complement of the positive halfspaces defined by each interval. It is a closed subset since each positive halfspace is open. The interiors of the translates of $D$ are all disjoint by the two previous propositions and the boundary components are identified pairwise, so $D$ is a fundamental domain for its orbit (Fig. \ref{fig:sp4halfspaces}). This orbit is in general hard to describe, but in some cases we can identify it precisely.

In the definition of generalized Schottky subgroups, we required that the model be a finite volume hyperbolization. This is an artificial requirement which made the analysis of maximal representations simpler. In what follows, we will assume that the model Schottky group acting on $\RP{1}$ is defined by intervals with disjoint closures, so it corresponds to an infinite volume hyperbolization. The advantage of using intervals with disjoint closures lies in the contraction property proven in \cite{Contraction} which we will exploit later.
%
%
%
%
\subsubsection{Intervals as symmetric spaces}
We will now describe how to identify an interval in $\Lag{2n}$ with the symmetric space associated with $\GL{n}$, endowing any interval with a Riemannian metric.\\
Let $P,Q \in \Lag{2n}$ be two transverse Lagrangians. As we saw earlier in Corollary \ref{Lem:PosTriplesTransverse}, all Lagrangians in the interval $(P,Q)$ have to be transverse to $Q$, so they are graphs of linear maps $f:P\to Q$. The isotropy condition on $f$ is given by
\[ \omega(v+f(v),v'+f(v')) = \omega(v,f(v')) + \omega(f(v),v') = 0 \quad \forall v,v' \in P. \]
Now we recall from our discussion of the Maslov index that we can associate the bilinear form
\begin{align*}
    \B{P}{Q}: P\oplus Q & \to \bR \\
    (v,w) & \mapsto \omega(v,\Inv{P}{Q}(w))
\end{align*}
to this splitting, and the index of its restriction to $\graph(f)$ is the Maslov index $\M(P,\graph(f),Q)$. We observe that this restriction is given by
\[ \B{P}{Q}(v+f(v),v'+f(v')) = \omega(v,f(v')) - \omega(f(v),v') = 2\omega(v,f(v')), \]
where the last equation follows from the isotropy condition on $f$. This bilinear form on $\graph(f)$ can also be seen as a symmetric bilinear form on $P$. Maximality of the Maslov index then translates to this form being positive definite.\\
Conversely, given a symmetric bilinear form $b$ on $P$, we obtain, for any $v'\in P$, a linear functional
\[ \left(v \mapsto \frac 12 b(v,v')\right) \in P^*. \]
Using the isomorphism
\begin{align*}
    Q & \to P^* \\
    w & \mapsto \omega(\cdot,w),
\end{align*}
we see that there is a unique vector $f(v') \in Q$ such that $b(v,v') = 2\omega(v,f(v')) \ \forall v\in P$. This uniquely defines a linear map $f:P\to Q$, and we have
\[ 2\left(\omega(v,f(v')) + \omega(f(v),v')\right) = b(v,v') - b(v',v) = 0, \]
so $\graph(f)$ is a Lagrangian. The Maslov index $\M(P,\graph(f),Q)$ is maximal if and only if $b$ is positive definite. This gives an identification of $(P,Q)$ with the space of positive definite symmetric bilinear forms on $P$, which is the symmetric space of $\mathrm{GL}(P)$.\\
The stabilizer in $\Sp{2n}$ of the pair $(P,Q)$ can be identified with $\mathrm{GL}(P)$ since any element $A\in\mathrm{GL}(P)$ uniquely extends to a symplectomorphism of $\bR^{2n}$ fixing $Q$: The linear forms $v \mapsto \omega(A(v),w)$ on $P$, for $w\in Q$, give rise to a unique automorphism $A^*:Q \to Q$ such that 
\[ \omega(A(v),w) = \omega(v,A^*(w)). \]
Then $A\oplus (B^*)^{-1}$ is the unique symplectic extension of $A$ fixing $Q$; we abuse notation slightly and denote it by $A$ as well. It acts on graphs $f:P\to Q$ by
\[ f \mapsto AfA^{-1}, \]
on bilinear forms on $P$ by
\[ (A\cdot b)(v,v') = b(A^{-1}v,A^{-1}v'), \]
and the identification of graphs and bilinear forms is equivariant with respect to these actions. In particular, $\mathrm{Stab}_{\Sp{2n}}(P,Q)$ is identified with the isometry group of the symmetric space $(P,Q)$.

\subsubsection{The Riemannian distance on intervals}
There is a simple formula for the Riemannian distance between two points in the interval $(P,Q)$:
\begin{Def}
Let $f,g$ be linear maps from $P$ to $Q$ whose graphs are elements of $(P,Q)$. Let $\lambda_1,\dots,\lambda_n$ be the eigenvalues of the automorphism $fg^{-1}$. Then, define
\[d_{PQ}(f,g) = \sqrt{\sum_{i=1}^n \log(\lambda_i)^2}.\]
\end{Def}
The following useful proposition is proved in \cite{Contraction}.
\begin{Prop}
\label{prop:contr}
Let $T\in\Sp{2n}$ such that $\overline{T(P,Q)} \subset (P,Q)$. Then, $T$ is a Lipschitz contraction for the distance $d_{PQ}$.
\end{Prop}
\begin{Cor}
\label{cor:contr}
Let $T\in \Sp{2n}$ such that $\overline{T(P,Q)} \subset (R,S)$. Then, for any $X,Y\in(P,Q)$,
\[d_{RS}(TX,TY) \leq C d_{PQ}(X,Y)\]
for some constant $0<C<1$.
\end{Cor}
Now we are ready to prove the main lemma that lets us describe the domains of discontinuity. Let $\rho:\Gamma\rightarrow\Sp{2n}$ define a generalized Schottky group in $\Sp{2n}$. We assume that the model $\Gamma$ is defined by intervals with distinct endpoints, so that the intervals in $\Lag{2n}$ have disjoint closures.
\begin{Lem} \label{lem:contr}
Let $\gamma \in \rho(\Gamma)$ be a word of reduced length $\ell$ in the generators $T_i$ and their inverses, with first letter $T$ and last letter $S$. We denote their attracting and repelling intervals by $I^\pm$ and $J^\pm$. Then, for any Schottky interval $K \neq J^-$ and $X,Y\in K$, we have
\[d_{I^+}(\gamma(X),\gamma(Y)) < C^{\ell} d_{K}(X,Y)\]
for some $0<C<1$ depending only on the set of generators.
\begin{proof}
Since a generator $T_k$ maps the interval $-I_k^-$ into $I_k^+$, we can consider it as a map from any Schottky interval $L\neq I_k^-$ into $I_k^+$. All of these maps are Lipschitz contractions by Corollary \ref{cor:contr}.\\
Now let $C$ be the maximum Lipschitz constant of all such maps, for $1\leq k\leq 2g$. We have
\[d_{J^+}(SX,SY) < C d_{K}(X,Y).\]
Composing contractions, we obtain
\[d_{I^+}(\gamma(X),\gamma(Y)) < C^\ell d_{K}(X,Y).\]
\end{proof}
\end{Lem}
\subsubsection{The domain of discontinuity}
We now have the necessary ingredients to analyze the orbit $\rho(\Gamma)\cdot D$ of the fundamental domain $D\subset\RP{2n-1}$ we defined in section \ref{sec:pos_halfspaces}. Using Lemma \ref{lem:contr}, we first define a map from the boundary of $\Gamma$, which is a Cantor set, into $\Lag{2n}$.
\begin{Prop}
Let $L\in \bigcap_{j=1}^{g}(-I_j^+\cap -I_j^-)$ where $I_j^\pm$ are the defining intervals of the generalized Schottky group. The evaluation map $\eta_0(\gamma) = \rho(\gamma)(L)$ induces an injective map $\eta:\partial \Gamma \rightarrow \Lag{2n}$ independent of the choice of $L$. Moreover, $\eta$ is continuous and increasing.
\begin{proof}
Let $x\in\partial\Gamma$ be a boundary point. Then $x$ corresponds to a unique infinite sequence in the generators $T_i$ and their inverses, where this sequence is reduced in the sense that no letter is followed by its inverse. We denote by $x^{(k)}\in\Gamma$ the word consisting of the first $k$ letters of $x$. Then the map $\eta$ will be defined by taking the limit
\[ \eta(x) = \lim\limits_{k\to\infty} \rho(x^{(k)})(L). \]
Let us first check that this limit does in fact exist. Recall that we introduced $k$-th order intervals and a bijection between words of length $k$ and $k$-th order intervals in section \ref{Sec:pingpong}. By the specific choice of $L$, its image $\rho(x^{(k)})(L)$ has to lie in the interval $I_{x^{(k)}}$ corresponding to the word $x^{(k)}$. Since the first $k$ letters of any word $x^{(m)}, \ m>k$ agree with $x^{(k)}$, the intervals $I_{x^{(k)}}$ form a nested sequence. Now we want to make use of the contraction property from the previous subsection. We first observe that since our model uses Schottky intervals with disjoint closures, second order intervals are relatively compact subsets of first order intervals. Let $I^{(2)} \subset I^{(1)}$ be such a configuration. Since the number of second order intervals is finite, there is a uniform bound $M$ such that
\[ \mathrm{diam}_{I^{(1)}}(I^{(2)}) < M, \]
where we used the metric on the symmetric space $I^{(1)}$. Then, denoting the first letter of $x$ by $T$, Lemma \ref{lem:contr} tells us that
\[ \mathrm{diam}_{I_T}(I_{x_k}) < MC^{k-2}. \]
This contracting sequence of nested subsets of the symmetric space $I_T$ thus has a unique limit, and $\eta$ is well-defined. By the same argument, we see that this limit does not depend on the choice of $L$.

We now show continuity of $\eta$. Let $y_n \to x$ be a sequence in $\partial\Gamma$ converging to $x$. This implies that for any $N \in \bN$, we can find $n_0$ such that for all $n\geq n_0$, the first $N$ letters of $y_n$ and $x$ agree. In this situation, $\eta(y_n)$ and $\eta(x)$ lie in the same interval $I_{x^{(N)}}$ and so we conclude, if the first letter of $x$ is $T$, that \[d_{I_T}(\eta(x),\eta(y_n)) < M C^{N-2}.\]

Finally, we prove positivity in a similar way to theorem \ref{Thm:bdry_map_gen_Schottky}. For any $x,y,z\in \partial\Gamma$ such that $\cycle{x}{y}{z}$ (where we use the natural embedding of $\partial \Gamma$ in $S^1$ to get the cyclic order) we can find a large enough $K$ so that $I_{x_K}$,$I_{y_K}$ and $I_{z_K}$ have disjoint closures. But since the cyclic relations on $k$-th order intervals are the same in $S^1$ as in $\Lag{2n}$, for any $P \in I_{x_K},Q\in I_{y_K},R\in I_{z_K}$ we have $\cycle{P}{Q}{R}$. In particular, $\cycle{\eta(x)}{\eta(y)}{\eta(z)}$.
\end{proof}
\end{Prop}
\begin{Rem}
The map $\eta$ that we define is related to the map $\xi$ of theorem \ref{Thm:bdry_map_gen_Schottky}. In this case, the endpoints of $k$-th order intervals are not dense in $S^1$, so we cannot get a map on the whole circle. However, because the intervals have disjoint closures, we get continuity on both sides rather than just left-continuity.
\end{Rem}
The next lemma relates the construction of the limit map $\eta$ with the positive halfspaces that intervals define in $\RP{2n-1}$.
\begin{Lem} \label{lem:conv_intervals}
Let $L_1^k,L_2^k$ be sequences of Lagrangians such that $L_1^k \rightarrow L$ and $L_2^k \rightarrow L$ with $\cycle{L_1^k}{L}{L_2^k}$ for all $k$. Then,
\[\bigcap_{k=1}^\infty \overline{\PH{L_1^k}{L_2^k}} = \bigcap_{k=1}^\infty \PH{L_1^k}{L_2^k} = \mathbb{P}(L).\]
\begin{proof}
Assume $\B{L_1^k}{L_2^k}(v,v)\geq 0$ for all $k$. Then we can find $v_k\xrightarrow{k\to\infty}v$ such that $\B{L_1^k}{L_2^k}(v_k,v_k)> 0$ for all $k$. Now, by proposition \ref{Prop:proj}, $v_k$ can be completed to a Lagrangian $L^k$ with $\M(L_1^k,L^k,L_2^k)= n$, so $L^k \subset (L_1^k,L_2^k)$ for all $k$, which implies $L^k \rightarrow L$, and so $v\in L$.
\end{proof}
\end{Lem}
Now we are ready to describe the orbit $\rho(\Gamma) D$.

The union of $D$ with the positive halfspaces defining the Schottky group is all of $\RP{2n-1}$, by definition of $D$. Denote by $\Gamma_\ell$ the set of words in $\Gamma$ of length up to $\ell$. Then, the union of $\rho(\Gamma_\ell) D$ with the projectivizations (positive halfspaces) of all $\ell$-th order intervals again covers all of $\RP{2n-1}$. Thus, when taking words of arbitrary length in $\Gamma$, these two pieces become respectively the full orbit $\rho(\Gamma) D$ and limits of nested positive halfspaces, which by lemma $\ref{lem:conv_intervals}$ collapse to the projectivization of a single Lagrangian. We conclude:
\begin{Thm}
\label{Thm:dod}
The orbit $\rho(\Gamma) D$ is the complement of a Cantor set of projectivized Lagrangian $n$-planes in $\RP{2n-1}$. This Cantor set is exactly the projectivization of the increasing set of Lagrangians defined by the boundary map $\eta$.
\end{Thm}
\begin{Rem}
The symplectic structure on $\bR^{2n}$ induces a contact structure on $\RP{2n-1}$ preserved by the symplectic group. The projectivizations of Lagrangian subspaces correspond to \emph{Legendrian} $(n-1)$-dimensional planes in $\RP{2n-1}$.
\end{Rem}

\bibliography{bibliography.bib}
\bibliographystyle{alpha}
\end{document}